\documentclass[reqno,11pt]{amsart}
\usepackage{amscd,amssymb}

\theoremstyle{plain}
\newtheorem{Thm}[subsection]{Theorem}
\newtheorem{Cor}[subsection]{Corollary}
\newtheorem{Lem}[subsection]{Lemma}
\newtheorem{Prop}[subsection]{Proposition}
\newtheorem{Conj}[subsection]{Conjecture}

\theoremstyle{definition}
\newtheorem{Def}[subsection]{Definition}

\theoremstyle{remark}

\newtheorem{Rem}[subsection]{Remark}

\newtheorem{conj}{Conjecture}

\numberwithin{equation}{section}
\renewcommand{\rm}{\normalshape}

\newif\ifShowLabels
\ShowLabelstrue
\newdimen\theight
\def\TeXref#1{%
    \leavevmode\vadjust{\setbox0=\hbox{{\tt
        \quad\quad  {\small \rm #1}}}%
    \theight=\ht0
    \advance\theight by \lineskip
    \kern -\theight \vbox to
    \theight{\rightline{\rlap{\box0}}%
    \vss}%
    }}%

\ShowLabelsfalse

\renewcommand{\sec}[2]{\section{#2}\label{S:#1}%
    \ifShowLabels \TeXref{{S:#1}} \fi}
\newcommand{\ssec}[2]{\subsection{#2}\label{SS:#1}%
    \ifShowLabels \TeXref{{SS:#1}} \fi}

\newcommand{\refs}[1]{Section ~\ref{S:#1}}

\newcommand{\reft}[1]{Theorem ~\ref{T:#1}}
\newcommand{\refl}[1]{Lemma ~\ref{L:#1}}

\newcommand{\refe}[1]{\eqref{E:#1}}
\newcommand{\refco}[1]{Conjecture ~\ref{Co:#1}}

\newenvironment{thm}[1]%
    { \begin{Thm} \label{T:#1}  \ifShowLabels \TeXref{T:#1} \fi }%
    { \end{Thm} }

\renewcommand{\th}[1]{\begin{thm}{#1} \sl }
\renewcommand{\eth}{\end{thm} }

\newenvironment{lemma}[1]%
    { \begin{Lem} \label{L:#1}  \ifShowLabels \TeXref{L:#1} \fi }%
    { \end{Lem} }
\newcommand{\lem}[1]{\begin{lemma}{#1} \sl}
\newcommand{\elem}{\end{lemma}}

\newenvironment{propos}[1]%
    { \begin{Prop} \label{P:#1}  \ifShowLabels \TeXref{P:#1} \fi }%
    { \end{Prop} }
\newcommand{\prop}[1]{\begin{propos}{#1}\sl }
\newcommand{\eprop}{\end{propos}}

\newenvironment{corol}[1]%
    { \begin{Cor} \label{C:#1}  \ifShowLabels \TeXref{C:#1} \fi }%
    { \end{Cor} }
\newcommand{\cor}[1]{\begin{corol}{#1} \sl }
\newcommand{\ecor}{\end{corol}}

\newenvironment{defeni}[1]%
    { \begin{Def} \label{D:#1}  \ifShowLabels \TeXref{D:#1} \fi }%
    { \end{Def} }
\newcommand{\defe}[1]{\begin{defeni}{#1} \sl }
\newcommand{\edefe}{\end{defeni}}

\newenvironment{remark}[1]%
    { \begin{Rem} \label{R:#1}  \ifShowLabels \TeXref{R:#1} \fi }%
    { \end{Rem} }
\newcommand{\rem}[1]{\begin{remark}{#1}}
\newcommand{\erem}{\end{remark}}

\newenvironment{conjec}[1]%
    { \begin{Conj} \label{Co:#1}  \ifShowLabels \TeXref{Co:#1} \fi }%
    { \end{Conj} }
\renewcommand{\conj}[1]{\begin{conjec}{#1} \sl }
\newcommand{\econj}{\end{conjec}}

\newcommand{\eq}[1]%
    { \ifShowLabels \TeXref{E:#1} \fi
       \begin{equation} \label{E:#1} }
\newcommand{\eeq}{ \end{equation} }

\newcommand{\prf}{ \begin{proof} }
\newcommand{\epr}{ \end{proof} }

\newcommand\alp{\alpha}		
		
\newcommand\gam{\gamma}		\newcommand\Gam{\Gamma}
\newcommand\del{\delta}		\newcommand\Del{\Delta}

\newcommand\lam{\lambda}		\newcommand\Lam{\Lambda}

\newcommand\calB{{\mathcal{B}}}

\newcommand\calF{{\mathcal{F}}}

\newcommand\calK{{\mathcal{K}}}
\newcommand\calL{{\mathcal{L}}}
\newcommand\calM{{\mathcal{M}}}

\newcommand\calO{{\mathcal{O}}}
\newcommand\calP{{\mathcal{P}}}

\newcommand\calS{{\mathcal{S}}}

\newcommand\calU{{\mathcal{U}}}

\newcommand\calW{{\mathcal{W}}}

\newcommand\QQ{\mathbb{Q}}

\newcommand\PP{\mathbb{P}}
\renewcommand\AA{\mathbb{A}}

\newcommand\GG{\mathbb{G}}

\newcommand\ZZ{\mathbb{Z}}

\newcommand\CC{\mathbb{C}}

	\newcommand\grg{{\mathfrak{g}}}

	\newcommand\grn{{\mathfrak{n}}}
	
	\newcommand\grp{{\mathfrak{p}}}

\newcommand\sdp{\times \hskip -0.3em {\raise 0.3ex
\hbox{$\scriptscriptstyle |$}}} 

\newcommand\Gr{\operatorname{Gr}}

\newcommand\Spec{\operatorname{Spec}}

\newcommand\Sym{\operatorname{Sym}}

\newcommand\Tr{\operatorname{Tr}}

\newcommand\os{{\overline{s}}}

\newcommand\oX{{\overline{X}}}

\newcommand\opi{{\overline{\pi}}}

\newcommand\hatB{{\widehat{B}}}

\newcommand\hatG{{\widehat{G}}}

\newcommand\hatU{{\widehat{U}}}

\newcommand\tils{{\widetilde{s}}}

\newcommand\x{\times}
\newcommand\ten{\otimes}

\newcommand{\ra}{\rangle}
\newcommand{\la}{\langle}

\newcommand\qlb{{\overline \QQ}_l}

\renewcommand\Spec{\operatorname{Spec}}

\newcommand{\kk}{\textsf{k}}
\newcommand\QM{\mathcal{Q M}}
\newcommand\ocalF{\overline{\calF}}
\newcommand\Fr{\operatorname{Fr}}
\newcommand\aff{\operatorname{aff}}
\newcommand\IC{\operatorname{IC}}
\newcommand\rank{\operatorname{rank}}
\begin{document}
Dedicated to S.~Patterson on the occasion of his 60th birthday

\bigskip

\author{Alexander Braverman, Michael Finkelberg and David Kazhdan}
\title{Affine Gindikin-Karpelevich formula via Uhlenbeck spaces}

\address{ A.B.: Department of Mathematics, Brown University,
151 Thayer St., Providence RI
02912, USA;
{\tt braval@math.brown.edu}}
\address{M.F.: IMU, IITP and State University Higher School of Economics,
Department of Mathematics, 20 Myasnitskaya st, Moscow 101000 Russia;
{\tt fnklberg@gmail.com}}
\address{D.K.: Einstein Institute of Mathematics, Edmond J. Safra Campus, Givat Ram The Hebrew
University of Jerusalem Jerusalem, 91904, Israel; {\tt kazhdan@math.huji.ac.il}}
\begin{abstract}We prove a version of the Gindikin-Karpelevich formula for untwisted affine Kac-Moody
groups over a local field of positive characteristic. The proof is geometric and it is based on the results
of \cite{bfg} about intersection cohomology of certain Uhlenbeck-type moduli spaces
(in fact, our proof is conditioned upon the assumption
that the results of \cite{bfg} are valid in positive characteristic; we believe that generalizing \cite{bfg} to the case of positive characteristic should be essentially
straightforward but we have not checked the details). In particular, we give a geometric
explanation of certain combinatorial differences
between finite-dimensional and affine case (observed earlier by Macdonald and Cherednik), which here manifest themselves
by the fact that the affine Gindikin-Karpelevich formula has an additional term compared to the finite-dimensional
case.  Very roughy speaking,
that additional term is related to the fact that
the loop group of an affine Kac-Moody group (which should be thought
of as some kind of ``double loop group") does not behave well from algebro-geometric point of view; however it has
a better behaved version which has something to do with algebraic surfaces.

A uniform (i.e. valid for all local fields) and unconditional (but not geometric) proof of the affine Gindikin-Karpelevich
formula is going to appear in \cite{BKP}.
\end{abstract}
\maketitle
\sec{}{The problem}
\ssec{}{Classical Gindikin-Karpelevich formula}
Let $\calK$ be a non-archimedian local field with ring of integers $\calO$ and let $G$ be a split semi-simple
group over $\calO$. The classical Gindikin-Karpelevich formula describes explicitly how a certain intertwining
operator acts on the spherical vector in a principal series representation of $G(\calK)$.
\footnote{More precisely, the Gindikin-Karpelevich
 formula answers the analogous question for real groups; its analog for $p$-adic groups (usually also referred to as Gindikin-Karpelevich formula) is proved e.g. in Chapter 4 of \cite{Langlands}.} In more explicit terms
it can be formulated as follows.

Let us choose a Borel subgroup $B$ of $G$ and an opposite Borel subgroup
$B_-$; let $U,U_-$ be their unipotent radicals. In addition, let  $\Lam$ denote the coroot lattice of $G$, $R_+\subset \Lam$ -- the set
of positive coroots, $\Lam_+$ -- the subsemigroup of $\Lam$ generated by $R_+$.
Thus any $\gam\in\Lam_+$ can be written as $\sum a_i\alp_i$ where $\alp_i$ are the simple roots.
We shall denote by $|\gam|$ the sum of all the $a_i$.

Set now $\Gr_G=G(\calK)/G(\calO)$. Then it is known that $\calU(\calK)$-orbits on $\Gr$ are in one-to-one correspondence
with elements of $\Lam$ (this correspondence will be reviewed in \refs{maps});
for any $\mu\in\Lam$ we shall denote by $S^{\mu}$ the corresponding orbit.
The same thing is true for $U_-(\calK)$-orbits. For each $\gam\in\Lam$ we shall denote by
$T^{\gam}$ the corresponding orbit. It is well-known that
$T^{\gam}\cap S^{\mu}$ is non-empty iff $\mu-\gam\in\Lam_+$ and in that case the above intersection
is finite. The Gindikin-Karpelevich formula allows one to compute the number of points
in $T^{-\gam}\cap S^0$ for $\gam\in\Lam_+$ (it is easy to see that the above intersection
is naturally isomorphic to $T^{-\gam+\mu}\cap S^{\mu}$ for any $\mu\in\Lam$). The answer is most easily stated
in terms of the corresponding generating function:
\th{main-finite}(Gindikin-Karpelevich formula)
$$
\sum\limits_{\gam\in\Lam_+} \# (T^{-\gam}\cap S^0)q^{-|\gam|}e^{-\gam}=
\prod\limits_{\alp\in R_+}\frac{1-q^{-1}e^{-\alp}}{1-e^{-\alp}}.
$$
\eth

\ssec{}{Formulation of the problem in the general case}
Let now $G$ be a split symmetrizable Kac-Moody group functor in the sense of \cite{Tits} and
let $\grg$ be the corresponding Lie algebra. We also let $\hatG$ denote the corresponding "formal" version
of $G$ (cf. page 198 in \cite{Tits}). The notations $\Lam,\Lam_+,R_+,\Gr_G,S^{\mu},T^{\gam}$
make sense for $\hatG$  without any changes (cf. \refs{maps} for more detail).
\conj{finite-general}
For any $\gam\in\Lam_+$ the intersection $T^{-\gam}\cap S^0$ is finite.
\econj
This conjecture will be proved in \cite{BKP} when $G$ is of affine type.
In this paper we are going to prove the following result:
\th{finite}
Assume that $\calK=\kk((t))$ where $\kk$ is finite. Then \refco{finite-general} holds.
\eth

So now (at least when $\calK$ is as above) we can ask the following

\medskip
\noindent
{\bf Question:} Compute the generating function\footnote{The reason that we use the notation $I_{\grg}$ rather than $I_{G}$ is that it is clear that this generating
function depends only on $\grg$ and not on $G$.}
$$
I_{\grg}(q)=\sum\limits_{\gam\in\Lam_+}\# (T^{-\gam}\cap S^0)~q^{-|\gam|}e^{-\gam}.
$$
One possible motivation for the above question is as follows: when $G$ is finite-dimensional,
Langlands \cite{Langlands} has observed that the usual
Gindikin-Karpelevich formula (more precisely, some generalization of it) is responsible for the fact that
the constant term of Eisenstein series induced from a parabolic subgroup of $G$ is related to some automorphic
$L$-function.
Thus we expect that generalizing the Gindikin-Karpelevich formula to general Kac-Moody group will eventially
become useful for studying Eisenstein series for those groups. This will be pursued in further publications.

We don't know the answer for general $G$. In the case when $G$ is finite-dimensional the answer is given by
\reft{main-finite}. In this paper we are going to reprove that formula by geometric means and give a generalization
to the case when $G$ is untwisted affine.
\ssec{}{The affine case}
Let us now assume that $\grg=\grg'_{\aff}$ where $\grg'$ is a simple finite-dimensional Lie algebra.
The Dynkin diagram of $\grg$ has a canonical ("affine") vertex and we let $\grp$ be the corresponding
maximal parabolic subalgebra of $\grg$. Let $\grg^{\vee}$ denote the Langlands dual algebra and
let $\grp^{\vee}$ be the corresponding dual parabolic. We denote by $\grn(\grp^{\vee})$ its
(pro)nilpotent radical.

Let $(e,h,f)$ be a principal $sl(2)$-triple in $(\grg')^{\vee}$. Since the Levi subalgebra
of $\grp^{\vee}$ is $\CC\oplus\grg'\oplus\CC$ (where the first multiple is central in $\grg^{\vee}$ and the
second is responsible for the ``loop rotation"), this triple acts on $\grn(\grp^{\vee})$ and we let
$\calW=(\grn(\grp^{\vee}))^f$ (the centralizer of $f$ in $\grn(\grp^{\vee})$). We are going to regard $\calW$ as a complex (with zero
differential) and with grading coming from the action of $h$ (thus $\calW$ is negatively graded).
In addition $\calW$ is endowed with an action of $\GG_m$, coming from the loop rotation in $\grg^{\vee}$.
In the case when $\grg'$ is simply laced we have $(\grg')^{\vee}\simeq \grg'$ and
$\grn(\grp^{\vee})=t\cdot\grg'[t]$ (i.e. $\grg'$-valued polynomials, which vanish at $0$).
Hence $\calW=t\cdot(\grg')^f[t]$ and the above $\GG_m$-action just acts by rotating
 $t$.  Let $d_1,\cdots,d_r$ be the exponents
of $\grg'$ (here $r=\rank(\grg')$). Then $(\grg')^f$ has a basis $(x_1,\cdots,x_r)$ where each $x_i$ is placed
in the degree $-2d_i$. We let $\Fr$ act on $\calW$ by requiring that it acts by
$q^{i/2}$ on elements of degree $i$. Also for any $n\in\ZZ$ let $\calW(n)$ be the same graded vector
space but with Frobenius action multiplied by $q^{-n}$.

Consider now $\Sym^*(\calW)$. We can again consider it as a complex concentrated in degrees $\leq 0$ endowed
with an action of $\Fr$ and $\GG_m$. For each $n\in\ZZ$ we let
$\Sym^*(\calW)_n$ be the part of $\Sym^*(\calW)$ on which $\GG_m$ acts by the character
$z\mapsto z^n$. This is a finite-dimensional complex with zero differential, concentrated in degrees
$\leq 0$ and endowed with an action of $\Fr$.

We are now ready to formulate the main result.
Let $\del$ denote the minimal positive imaginary coroot of $\grg$.
Set
$$
\Del_{\calW}(z)=\sum\limits_{n=0}^{\infty}\Tr(\Fr,\Sym^*(\calW)_n)z^n.
$$
In particular, when $\grg'$ is simply laced we have
$$
\Del(z)=\prod\limits_{i=1}^r \prod\limits_{j=0}^{\infty}(1-q^{-d_i}z^j)^{-1}.
$$

\th{main-affine}(Affine Gindikin-Karpelevich formula)

\noindent
Assume that the results of
\cite{bfg} are valid over $\kk$ and let $\calK=\kk((t))$. Then
$$
I_{\grg}(q)=\frac{\Del_{\calW(1)}(e^{-\del})}{\Del_{\calW}(e^{-\del})}\prod\limits_{\alp\in R_+}\left(\frac{1-q^{-1}e^{-\alp}}{1-e^{-\alp}}\right)^{m_{\alp}}.
$$
Here $m_\alp$ denotes the multiplicity of the coroot $\alp$.
\eth

\noindent
{\bf Remark.} Although formally the paper \cite{bfg} is written under the assumption
that $\text{char}~ \kk=0$, we believe that adapting all the constructions of \cite{bfg}
to the case
$\text{char}~ \kk=p$ should be more or less straightforward. We plan to discuss it in a separate publication.

\medskip
\noindent
Let us make two remarks about the above formula: first, we see that it is very similar
to the finite-dimensional case (of course in that case $m_\alp=1$ for any $\alp$)
with the exception of a ``correction term" (which is equal to
$\frac{\Del_{\calW}(e^{-\del})}{\Del_{\calW(1)}(e^{-\del})}$). Roughly speaking this correction term has to do with
imaginary coroots of $\grg$.
The second remark is that the same correction term appeared in the work of Macdonald
\cite{Mac} from purely combinatorial point of view (cf. also \cite{CherMa} for a more detailed study). The
main purpose of this note is to explain how the term $\frac{\Del_{\calW}(e^{-\del})}{\Del_{\calW(1)}(e^{-\del})}$ appears naturally from
geometric point of view (very roughly speaking it is related to the fact that affine Kac-
Moody groups over a local field of positive characteristic can be studied using various
moduli spaces of bundles on an algebraic surface). The relation between the present work
and the constructions of \cite{CherMa} and \cite{Mac}  will be discussed in \cite{BKP}.
\ssec{}{Acknowledgments}We thank I.~Cherednik, P.~Etingof and M.~Patnaik for very helpful discussions.
A.~B. was partially supported by the NSF grant DMS-0901274.
M.~F. was partially supported by the RFBR grant 09-01-00242 and the Science
Foundation of the SU-HSE awards No.T3-62.0 and 10-09-0015.
D.~K. was partially supported by the BSF grant
037.8389.

\sec{maps}{Interpretation via maps from $\PP^1$ to $\calB$}

\ssec{}{Generalities on Kac-Moody groups}
In what follows all schemes will be considered over a field $\kk$ which at some point will be assumed to be finite.
Our main reference for Kac-Moody groups is \cite{Tits}.
Assume that we are given a symmetrizable Kac-Moody root data and we denote by $G$ (resp. $\hatG$) the corresponding
minimal (resp. formal) Kac-Moody group functor (cf. \cite{Tits}, page 198); we have the natural
embedding $G\hookrightarrow \hatG$. We also let $W$ denote the corresponding Weyl group and we let
$\ell:W\to \ZZ_{\geq 0}$ be the corresponding length function.

The group $G$ is endowed with closed subgroup functors $U\subset B, U_-\subset B_-$ such that
the quotients $B/U$ and $B_-/U_-$ are naturally isomorphic to the Cartan group $H$ of $G$; also
$H$ is isomorphic to the intersection $B\cap B_-$.
Moreover, both $U_-$ and $B_-$ are still closed as subgroup functors of $\hatG$. On the other hand,
$B$ and $U$ are not closed in $\hatG$ and we denote by $\hatB$ and $\hatU$ their closures.

The quotient $G/B$ has a natural structure of an ind-scheme which is ind-proper; the same is true
for the quotient $\hatG/\hatB$ and the natural
map $G/B\to \hatG/\hatB$ is an isomorphism. This quotient is often called {\em the thin flag variety of $G$}.
Similarly, one can consider the quotient $\calB=\hatG/B_-$; it is called {\em the thick flag variety of $G$}
or {\em Kashiwara flag scheme}. As is suggested by the latter name, $\calB$ has a natural scheme structure.
The orbits of $B$ on $\calB$ are in one-to-one correspondence with the elements of the
Weyl group $W$; for each $w\in W$ we denote by $\calB_w$ the corresponding orbit.
The codimension of $\calB_w$ is $\ell(w)$; in particular, $\calB_e$ is open.
There is a unique $H$-invariant point $y_0\in\calB_e$. The complement to $\calB_e$ is a divisor in $\calB$ whose
components are in one-to-one correspondence with the simple roots of $G$.

In what follows $\Lam$ will denote the coroot lattice of $G$, $R_+\subset \Lam$ -- the set
of positive coroots, $\Lam_+$ -- the subsemigroup of $\Lam$ generated by $R_+$.
Thus $\gam\in\Lam_+$ can be written as $\sum a_i\alp_i$ where $\alp_i$ are the simple coroots.
We shall denote by $|\gam|$ the sum of all the $a_i$.

In what follows we shall assume that $G$ is "simply connected", which means that
$\Lam$ is equal to the full cocharacter lattice of $H$.
\ssec{}{Some further notations}
For any variety $X$ and any $\gam\in\Lam_+$ we shall denote by $\Sym^{\gam}X$ the variety parametrizing
all unordered collections $(x_1,\gam_1),...(x_n,\gam_n)$ where $x_j\in X, \gam_j\in\Lam_+$ such that
$\sum\gam_j=\gam$.

Assume that $\kk$ is finite and let $\calS$ be a complex of $\ell$-adic sheaves on a variety $X$ over
$\kk$.
We set
$$
{\chi}_{_\kk}(\calS)=\sum\limits_{i\in\ZZ}(-1)^i\Tr(\Fr, H^i(\oX,\calS)),
$$
where $\oX=X\underset{\Spec{\kk}}\x\Spec{\overline{\kk}}$.

We shall denote by $(\qlb)_X$ the constant sheaf with fiber $\qlb$. According to the Grothendieck-Lefschetz
fixed point formula we have
$$
{\chi}_{_\kk}((\qlb)_X)=\# X(\kk).
$$
\ssec{}{Semi-infinite orbits}As in the introduction we set $\calK=\kk((t))$, $\calO=\kk[[t]]$.
We let $\Gr=\hatG(\calK)/\hatG(\calO)$, which we are just going to consider as a set with no structure.
 Each $\lam\in \Lam$ is a homomorphism
$\GG_m\to H$; in particular, it defines a homomorphism $\calK^*\to H(\calK)$. We shall
denote the image of $t$ under the latter homomorphism by $t^{\lam}$. Abusing the notation, we shall
denote its image in $\Gr$ by the same symbol.
Set
$$
S^{\lam}=\hatU(\calK)\cdot t^{\lam}\subset \Gr;\quad T^{\lam}=U_-(\calK)\cdot t^{\lam}\subset \Gr.
$$
\lem{iwasawa} $\Gr$ is equal to the disjoint union of all the $S^{\lam}$.
\elem
\prf This follows from the
{\em Iwasawa decomposition} for $G$ of \cite{Gauss-Rous}; we include a different proof for completeness.
Since $\Lam\simeq \hatU(\calK)\backslash \hatB(\calK)/\hatB(\calO)$,
the statement of the lemma is equivalent to the assertion that
the natural map $\hatB(\calK)/\hatB(\calO)\to \hatG(\calK)/\hatG(\calO)$ is an isomorphism; in other words, we need
to show that $\hatB(\calK)$ acts transitively on $\Gr$. But this is equivalent to saying that
$\hatG(\calO)$ acts transitively on $\hatG(\calK)/\hatB(\calK)$, which means that the natural map
$\hatG(\calO)/\hatB(\calO)\to \hatG(\calK)/\hatB(\calK)$ is an isomorphism. However, the left hand side
is $(\hatG/\hatB)(\calO)$ and the right hand side is $(\hatG/\hatB)(\calK)$ and the assertion follows
from the fact that the ind-scheme $\hatG/\hatB$ satisfies the valuative criterion of properness.
\epr

\medskip
\noindent
The statement of the lemma is definitely false if we use $T^{\mu}$'s instead of $S^{\lam}$'s since
the scheme $\hatG/B_-$ does not satisfy the valuative ctiterion of properness.
Let us say that an element $g(t)\in \hatG(\calK)$ is good if its projection
to $\calB(\calK)=B_-(\calK)\backslash\hatG(\calK)$ comes from a point of $\calB(\calO)$. Since
$\calB(\calO)=B_-(\calO)\backslash\hatG(\calO)$, it follows that the set of good elements of
$\hatG(\calK)$ is just equal to $B_-(\calK)\cdot G(\calO)$, which immediately proves the following
result:
\lem{iwasawa-}
The preimage of $\bigcup\limits_{\gam\in \Lam}T^{\gam}$ in $\hatG(\calK)$ is equal to the set of good elements
of $\hatG(\calK)$.
\elem
\ssec{}{Spaces of maps}Recall that the Picard group of $\calB$ can be naturally identified with
$\Lam^{\vee}$ (the dual lattice to $\Lam$). Thus for any map $f:\PP^1\to\calB$ we can talk about the
degree of $f$ as an element $\gam\in\Lam$. The space of such maps is non-empty iff $\gam\in\Lam_+$.
We say that a map $f:\PP^1\to\calB$ is {\em based} if $f(\infty)=y_0$.
Let $\calM^{\gam}$ be the space of based maps $f:\PP^1\to \calB$ of degree $\gam$. It is shown in the Appendix
to \cite{bfg} that this is a smooth scheme of finite type over $\kk$ of dimension $2|\gam|$.
We have a natural (``factorization") map $\pi^{\gam}:\calM^{\gam}\to\Sym^{\gam}\AA^1$ which is related to how
the image of a map $\PP^1\to\calB$ intersects the complement to $\calB_e$. In particular, if we set
$$
\calF^{\gam}=(\pi^{\gam})^{-1}(\gam\cdot 0),
$$
then $\calF^{\gam}$ consists of all the based maps $f:\PP^1\to\calB$ of degree $\gam$ such that
$f(x)\in\calB_e$ for any $x\neq 0$.
\th{fgamma}
There is a natural identification $\calF^{\gam}(\kk)\simeq T^{-\gam}\cap S^0$.
\eth
Since $\calF^{\gam}$ is a scheme of finite type over $\kk$, it follows that $\calF^{\gam}(\kk)$ is finite and thus
\reft{fgamma} implies \reft{finite}.

The proof of \reft{fgamma} is essentially a repetition of a similar proof in the finite-dimensional case, which we include here for
completeness.
\prf
First of all, let us construct an embedding of the union of all the $\calF^{\gam}(\kk)$ into $S^0=\hatU(\calK)/\hatU(\calO)$.
Indeed an element of $\bigcup\limits_{\gam\in\Lam_+}\calF^{\gam}$ is uniquely determined by its restriction
to $\GG_m\subset \PP^1$; this restriction is a map
$f:\GG_m\to\calB_e$ such that $\lim_{x\to\infty} f(x)=y_0$. We may identify $\calB_e$ with
$\hatU$ (by acting on $y_0$). Thus we get
\eq{identification}
\bigcup\limits_{\gam\in\Lam_+}\calF^{\gam}\subset \{u:\PP^1\backslash\{0\}\to\hatU|\ u(\infty)=e\}.
\end{equation}
We have a natural map from the set of $\kk$-points of the right hand side of \refe{identification}
to $\hatU(\calK)$; this map sends every $u$ as above to its restriction to the formal punctured neighbourhood of
$0$. We claim that after projecting $\hatU(\calK)$ to $S^0=\hatU(\calK)/\hatU(\calO)$, this map becomes an isomorphism.
Recall that $\hatU$ is a group-scheme, which can be written as a projective limit of finite-dimensional unipotent group-schemes
$U_i$; moreover, each $U_i$ has a filtration by normal subgroups with successive quotients isomorphic to
$\GG_a$. Hence it is enough to prove that the above map is an isomorphism when $U=\GG_a$.
In this case we just need to check that any element of the quotient $\kk((t))/\kk[[t]]$ has
unique lift to a polynomial $u(t)\in \kk[t,t^{-1}]$ such that $u(\infty)=0$, which is obvious.

Now \refl{iwasawa-} implies that a map $u(t)$ as above extends to a map $\PP^1\to\calB$ if and only
if the corresponding element of $S^0$ lies in the intersection with some $T^{-\gam}$.

It remains to show that $\calF^{\gam}(\kk)$ is exactly equal to $S^0\cap T^{-\gam}$ as a subset of $S^0$.
 Let $\Lam^{\vee}$ be the weight lattice of $G$ and let
$\Lam_+^{\vee}$ denote the set of dominant weights of $G$. For each $\lam^{\vee}\in\Lam_+^{\vee}$
we can consider the Weyl module $L(\lam^{\vee})$, defined over $\ZZ$;  in particular, $L(\lam^{\vee})(\calK)$ and
$L(\lam^{\vee})(\calO)$ make sense. By the definition $L(\lam^{\vee})$ is the module
of global sections of a line bundle $\calL(\lam^{\vee})$ on $\calB$. Moreover we have a weight decomposition
$$
L(\lam^{\vee})=\bigoplus\limits_{\mu^{\vee}\in\Lam^{\vee}} L(\lam^{\vee})_{\mu^{\vee}}
$$
where each $L(\lam^{\vee})_{\mu^{\vee}}$ is a finitely generated free $\ZZ$-module and
$L(\lam^{\vee})_{\lam^{\vee}}:=l_{\lam^{\vee}}$ has rank one. Geometrically, $l_{\lam^{\vee}}$ is the
fiber of $\calL(\lam^{\vee})$ at $y_0$ and the corresponding projection map from $L(\lam^{\vee})=\Gam(\calB,\calL(\lam^{\vee}))$
to $l_{\lam^{\vee}}$ is the restriction to $y_0$.

Let $\eta_{\lam^{\vee}}$ denote the projection of $L(\lam^{\vee})$ to $l_{\lam^{\vee}}$. This map is
$U_-$-equivariant (where $U_-$ acts trivially on $l_{\lam^{\vee}}$).
\lem{}
The projection of a good element $g\in G(\calK)$ lies in $T^{\nu}$ (for some $\nu\in\Lam$) if and only if
for any $\lam^{\vee}\in\Lam^{\vee}$ we have:
\eq{condition}
\eta_{\lam^{\vee}}(g(L(\lam^{\vee})(\calO)))\subset t^{\la\nu,\lam^{\vee}\ra}l_{\lam^{\vee}}(\calO);\qquad
\eta_{\lam^{\vee}}(g(L(\lam^{\vee})(\calO)))\not\subset t^{\la\nu,\lam^{\vee}\ra-1}l_{\lam^{\vee}}(\calO).
\end{equation}
\elem
\prf
First of all, we claim that if the projection of $g$ lies in $T^{\nu}$ then the above condition is satisfied.
Indeed, it is clearly satisfied by $t^{\nu}$; moreover \refe{condition} is clearly invariant under left
multiplication by $U_-(\calK)$ and under right multiplication by $G(\calO)$.
Hence any $g\in U_-(\calK)\cdot t^{\nu}\cdot G(\calO)$ satisfies \refe{condition}.

On the other hand, assume that a good element $g\in G(\calK)$ satisfies \refe{condition}. Since $g$ lies in $U_-(\calK)\cdot t^{\nu'}\cdot G(\calO)$ for some $\nu'$, it follows that $g$
satisfies \refe{condition} when $\nu$ is replaced by $\nu'$. However, it is clear that this is possible only if $\nu=\nu'$.
\epr

It is clear that in \refe{condition} one can replace $g(L(\lam^{\vee})(\calO))$ with
$g(L(\lam^{\vee})(\kk))$ (since the latter generates the former as an $\calO$-module).

Let now $f$ be an element of $\calF^{\gam}$. Then $f^*\calL(\lam^{\vee})$ is isomorphic to the line bundle
$\calL(\la \gam,\lam^{\vee}\ra)$ on $\PP^1$. On the other hand, the bundle $\calL(\lam^{\vee})$ is trivialized
on $\calB_e$ by means of the action of $U$; more precisely, the restriction of $\calL(\lam^{\vee})$ is canonically
identified with the trivial bundle with fiber $l_{\lam^{\vee}}$. Let now $s\in L(\lam^{\vee})(\kk)$; we are going
to think of it as a section of $L(\lam^{\vee})$ on $\calB$. In particular, it gives rise to a function
$\tils:\calB_e\to l_{\lam^{\vee}}$. Let also $u(t)$ be the element of $U(\calK)$, corresponding to $f$. Then $\eta_{\lam^{\vee}}(u(t)(s))$ can be described as follows: we consider the composition $\tils\circ f$
and restrict it to the formal neighbourhood of $0\in \PP^1$ (we get an element of $l_{\lam^{\vee}}(\calK)$).
On the other hand, since $f\in\calF^{\gam}$, it follows that $f^*\calL(\lam^{\vee})$
is trivialized away from $0$ and any section of it can be thought of as a function $\PP^1\backslash \{0\}$ with
pole of order $\leq \la \gam,\lam^{\vee}\ra$ at $0$. Hence $\tils\circ f$ has pole of order $\leq \la \gam,\lam^{\vee}\ra$
at $0$.

To finish the proof it is enough to show that for some $s$ the function $\tils\circ f$ has pole of order exactly $\la \gam,\lam^{\vee}\ra$ at $0$ (indeed if $f\in T^{-\gam'}$ for some $\gam'\in \Lam$, then by \refe{condition}
$\tils\circ f$ has pole of order $\leq \la \gam',\lam^{\vee}\ra$ at $0$ and for some $s$, it has pole of order exactly
$ \la \gam',\lam^{\vee}\ra$ which implies that $\gam=\gam'$). To prove this, let us note that since $\calL(\lam^{\vee})$
is generated by global sections, the line bundle $f^*\calL(\lam^{\vee})$ is generated by sections of the form
$f^*s$, where $s$ is a global section of $\calL(\lam^{\vee})$. This implies that for any
$\os\in \Gam(\PP^1,f^*\calL(\lam^{\vee}))$ there exists a section $s\in\Gam(\calB,\calL(\lam^{\vee}))$ such that
the ratio $s/\os$ is a rational function on $\PP^1$, which is invertible at $0$. Taking $\os$ such that
its pole with respect to the above trivialization of $f^*\calL(\lam^{\vee})$ is exactly equal to
$ \la \gam',\lam^{\vee}\ra$ and taking $s$ as above, we see that the pole of $f^*s$ with respect to the above trivialization of $f^*\calL(\lam^{\vee})$ is exactly equal to
$ \la \gam',\lam^{\vee}\ra$.
\epr
\sec{finite}{Proof of \reft{main-finite} via quasi-maps}
\ssec{}{Quasi-maps}
We shall denote by $\QM^{\gam}$ the space of based {\em quasi-maps} $\PP^1\to \calB$.
According to \cite{ffkm}
we have the stratification
$$
\QM^{\gam}=\bigcup\limits_{\gam'\leq \gam}\calM^{\gam'}\x\Sym^{\gam-\gam'}\AA^1.
$$
The factorization morphism $\pi^{\gam}_0$ extends to the similar morphism
$\opi^\gam:\QM^{\gam}\to \Sym^{\gam}$ and we set $\ocalF^{\gam}=(\opi^{\gam})^{-1}(0)$.
Thus we have
\eq{decomposition}
\ocalF^{\gam}=\bigcup\limits_{\gam'\leq\gam}\calF^{\gam'}.
\end{equation}
There is a natural section $i^{\gam}:\Sym^{\gam}\AA^1\to \QM^{\gam}$. According to \cite{ffkm} we have
\th{ic-fink}
\begin{enumerate}
\item
The restriction of $\IC_{\QM^{\gam}}$ to $\calF^{\gam'}$ is isomorphic
to $(\qlb)_{\calF^{\gam'}}[2](1)^{\ten|\gam'|}\ten \Sym^*(\grn_+^{\vee}[2](1))_{\gam-\gam'}$.
\item
There exists a $\GG_m$-action on $\QM^{\gam}$ which contracts it to the image of $i^{\gam}$. In particular,
it contracts $\ocalF^{\gam}$ to one point (corresponding to $\gam'=0$ in \refe{decomposition}).
\item
Let $s_{\gam}$ denote the embedding of $\gam\cdot 0$ into $\Sym^{\gam}\AA^1$. Then
$$
s_{\gam}^* i_{\gam}^!\IC_{\QM^{\gam}}=\Sym^*(\grn_+)_{\gam}
$$
(here the right hand is a vector space concentrated in cohomological degree 0 and with trivial action
of $\Fr$).
\end{enumerate}
\eth
The assertion 2) implies that
$\pi^{\gam}_!\IC_{\QM^{\gam}}=i_{\gam}^!\IC_{\QM^{\gam}}$ and hence
$$
H^*_c(\ocalF,\IC_{\QM^{\gam}}|_{\ocalF^{\gam}})=s_{\gam}^*\pi^{\gam}_!\IC_{\QM^{\gam}}=s_{\gam}^*i_{\gam}^!\IC_{\QM^{\gam}}=
\Sym^*(\grn_+)_{\gam}.
$$
Thus, setting, $\calS^{\gam}=\IC_{\QM^{\gam}}|_{\ocalF^{\gam}}$ we get
\eq{one}
\sum\limits_{\gam\in\Lam_+}{\chi}_{_\kk}(\calS^{\gam})e^{-\gam}=\prod\limits_{\alp\in R_+}\frac{1}{1-e^{-\alp}}.
\end{equation}
On the other hand, according to 1) we have
$$
{\chi}_{_\kk}(\calS^{\gam})=\sum\limits_{\gam'\leq\gam}(\#\calF^{\gam'})q^{-|\gam'|}
\Tr(\Fr,\Sym^*(\grn_+^{\vee}[2](1))_{\gam-\gam'})
$$
which implies that
\eq{two}
\sum\limits_{\gam\in\Lam_+}{\chi}_{_\kk}(\calS^{\gam})e^{-\gam}=
\frac{\sum\limits_{\gam\in\Lam_+}\#\calF^{\gam}(\kk)q^{-|\gam|}e^{-\gam}}
{\prod\limits_{\alp\in R_+}{1-q^{-1}e^{-\alp}}}=
\frac{I_{\grg}(q)}{\prod\limits_{\alp\in R_+}{1-q^{-1}e^{-\alp}}}.
\end{equation}
Hence
$$
I_{\grg}(q)=\prod\limits_{\alp\in R_+}\frac{1-q^{-1}e^{-\alp}}{1-e^{-\alp}}.
$$
\sec{affine}{Proof of \reft{main-affine}}
\ssec{}{Flag Uhlenbeck spaces}
We now assume that $G=(G')_{\aff}$ where $G'$ is some semi-simple simply connected group.
We want to follow the pattern of \refs{finite}.
Let $\gam\in\Lam_+$. As is discussed
in \cite{bfg} the corresponding space of quasi-maps behaves badly when $G$ is replaced by $G_{\aff}$.
However, in this case one can use the corresponding {\em flag Uhlenbeck space} $\calU^{\gam}$.
In fact, as was mentioned in the Introduction,  in \cite{bfg} only the case of $\kk$ of characteristic 0
is considered. In what follows we are going to assume that the results of
{\em loc. cit.} are valid also in positive characteristic.

The flag Uhlenbeck space $\calU^{\gam}$ has properties similar to the space of quasi-maps $\QM^{\gam}$ considered above in the previous Section.
Namely we have:

\medskip
a) $\calU^{\gam}$ is an affine variety of dimension $2|\gam|$, which contains $\calM^{\gam}$ as a dense open subset.

b) There is a factorization map $\pi^{\gam}:\calU^{\gam}\to \Sym^{\gam}\AA^1$; it has a section
$i_{\gam}:\Sym^{\gam}\AA^1\to \calU^{\gam}$.

c)  $\calU^{\gam}$ is endowed with a $\GG_m$-action which contracts $\calU^{\gam}$ to the image
of $i_{\gam}$.

\medskip
\noindent
These properties are identical to the corresponding properties of $\QM^{\gam}$ from the previous Section.
The next (stratification) property, however, is different (and it is in fact responsible for
the additional term in \reft{main-affine}). Namely, let $\del$ denote the minimal positive imaginary
coroot of $G'_{\aff}$. Then we have

\medskip
d) There exists a stratification
\eq{uhl-strat}
\calU^{\gam}=\bigcup\limits_{\gam'\in\Lam_+, n\in\ZZ,\gam-\gam'-n\del\in\Lam_+ }
(\calM^{\gam-\gam'-n\del}\x\Sym^{\gam'}\AA^1)\x \Sym^n (\GG_m\x\AA^1).
\end{equation}

\medskip
\noindent
In particular, if we now set $\ocalF^{\gam}=(\pi^{\gam})^{-1}(\gam\cdot 0)$ we get
\eq{str-fgam}
\ocalF^{\gam}=(\bigcup\limits_{\gam'\in\Lam_+, n\in\ZZ,\gam-\gam'-n\del\in\Lam_+ }
\calF^{\gam'})\x\Sym^n(\GG_m).
\end{equation}
\ssec{}{Description of the IC-sheaf}In \cite{bfg} we describe the IC-sheaf of $\calU^{\gam}$.
To formulate the answer, we need to introduce some notation. Let
$\calP(n)$ denote the set of partitions of $n$. In other words, any $P\in \calP(n)$ is an unordered sequence
$n_1,\cdots,n_k\in \ZZ_{>0}$ such that $\sum n_i=n$. We set $|P|=k$. For a variety $X$ and any $P\in\calP(n)$ we denote
by $\Sym^P(X)$ the locally closed subset of $\Sym^n(X)$ consisting of all formal sums
$\sum n_i x_i$ where $x_i\in X$ and $x_i\neq x_j$ for $i\neq j$. The dimension of $\Sym^P(X)$ is $|P|\cdot\dim X$. Let also
$$
\Sym^*(W[2](1))_P=\bigotimes\limits_{i=1}^k\Sym^*(\calW[2](1))_{n_i}.
$$
\th{bfg-ic}
The restriction of $\IC_{\calU^{\gam}}$ to $\calM^{\gam-\gam'-n\del}\x\Sym^{\gam'}(\AA^1)\x\Sym^P(\GG_m\x\AA^1)$
is isomorphic to constant sheaf on that scheme tensored with
$$
\Sym^*(\grn_+)_{\gam'}\ten \Sym^*(\calW)_P[2|\gam-\gam'-n\del|](|\gam-\gam'-n\del|).
$$
\eth
\cor{bfg-cor}The restriction of $\IC_{\calU^{\gam}}$ to $\calF^{\gam-\gam'-n\del}\x\Sym^P(\GG_m)$ is isomorphic to
the constant sheaf tensored with
$$
\Sym^*(\grn_+^{\vee})_{\gam'}\ten \Sym^*(\calW)_P[2|\gam-\gam'-n\del|](|\gam-\gam'-n\del|).
$$
\ecor

Let now $\calS^{\gam}$ denote the restriction of $\IC_{\calU^{\gam}}$ to $\ocalF^{\gam}$.
Then, as in \refe{one} we get
\eq{one-uhl}
\sum\limits_{\gam\in\Lam_+}{\chi}_{_\kk}(\calS^{\gam})e^{-\gam}=\prod\limits_{\alp\in R_+}\frac{1}{(1-e^{-\alp})^{m_\alp}}.
\end{equation}
On the other hand, arguing as in \refe{two} we get
that
\eq{two-uhl}
\sum\limits_{\gam\in\Lam_+}{\chi}_{_\kk}(\calS^{\gam})e^{-\gam}=
A(q)\frac{I_{\grg}(q)}{\prod\limits_{\alp\in R_+}{(1-q^{-1}e^{-\alp})^{m_\alp}}}.
\end{equation}
where
$$
A(q)=\sum\limits_{n=0}^{\infty}\sum\limits_{P\in\calP(n)}\Tr(\Fr, H^*_c(\Sym^P(\GG_m),\qlb)\ten\Sym^*(\calW[2](1))_P)e^{-n\del}.
$$
This implies that
$$
I_{\grg}(q)=A(q)^{-1}\prod\limits_{\alp\in R_+}\left(\frac{1-q^{-1}e^{-\alp}}{1-e^{-\alp}}\right)^{m_\alp}.
$$
It remains to compute $A(q)$. However, it is clear that
\eq{bredd}
A(q)=\sum\limits_{n=0}^{\infty} \Sym^n(H^*_c(\GG_m)\ten \calW[2](1))e^{-n\del}=\frac{\Del_{\calW}(e^{-\del})}{\Del_{\calW(1)}(e^{-\del})}.
\end{equation}
This is true since $H^i_c(\GG_m)=0$ unless $i=1,2$, and we have
$$
H^1_c(\GG_m)=\qlb,\quad H^2_c(\GG_m)=\qlb(-1),
$$
and thus if we ignore the cohomological $\ZZ$-grading, but only remember the corresponding $\ZZ_2$-grading,
then we just have
$$
\Sym^*(H^*_c(\GG_m)\ten \calW[2](1))=\Sym^*(\calW)\ten \Lambda^*(\calW(1)),
$$
whose character is exactly the right hand side of \refe{bredd}.

\end{document}